\documentclass[12pt]{article}
\usepackage{amsfonts}
\usepackage{amsmath}
\usepackage{graphicx}

\input tcilatex
\input tcilatex
\begin{document}

\author{Ozlem Ersoy \thanks{%
Ozlem Ersoy, Telephone: +90 222 2393750, Fax: +90 222 2393578, E-mail:
ozersoy@ogu.edu.tr} and Idiris Dag \\
Mathematics-Computer Department,\\
Eski\c{s}ehir Osmangazi University, 26480, Eski\c{s}ehir, Turkey.$^{1}$}
\title{A Trigonometric Cubic B-spline Finite Element Method for Solving the
Nonlinear Coupled Burger Equation}
\maketitle

\begin{abstract}
The coupled Burgers equation is solved by way of the trigonometric B-spline
collocation method. The unknown of the coupled Burgers equation is
integrated in time by aid of the Crank-Nicolson method. Resulting
time-integrated coupled Burgers equation is discretized using the
trigonometric cubic B-spline collocation method. Fully-integrated couupled
Burgers equation which is a system of nonlinear algebraic equation is solved
with a variant of Thomas algorithm. The three model test problems are
studied to illustrate the accuracy of the suggested method.
\end{abstract}

\section{Introduction}

\qquad One of the flow equations, derived by Episov\cite{1995}, is the
coupled Burgers equation(CBE) describing sedimentation or evaluation of the
scaled volume concentration of two kinds particles in fluid suspensions or
colloids under the effect of gravity\cite{1998}. \ The CBE is known as the
simple case of the Navier-Stokes equation due to including the nonlinear
convection term and viscosity term

The following \ nonlinear partial differential equation describes the CBE%
\begin{equation}
\begin{array}{c}
\dfrac{\partial U}{\partial t}-\dfrac{\partial ^{2}U}{\partial x^{2}}+k_{1}U%
\dfrac{\partial U}{\partial x}+k_{2}(UV)_{x}=0 \\ 
\\ 
\dfrac{\partial V}{\partial t}-\dfrac{\partial ^{2}V}{\partial x^{2}}+k_{1}V%
\dfrac{\partial V}{\partial x}+k_{3}(UV)_{x}=0%
\end{array}
\label{CB1}
\end{equation}%
where $k_{1}$, $k_{2}$ and $k_{3}$ are reel constants and subscripts $x$ and 
$t$ denote differentiation, $x$ distance and $t$ time, is considered.
Boundary conditions%
\begin{equation}
\text{\ }%
\begin{array}{l}
U(a,t)=f_{1}(a,t),\text{ }U(b,t)=f_{2}(b,t) \\ 
V(a,t)=g_{1}(a,t),\text{ }V(b,t)=g_{2}(b,t),\text{ }t>0%
\end{array}
\label{CB2}
\end{equation}%
and initial conditions%
\begin{eqnarray}
U(x,0) &=&f(x)  \label{CB3} \\
V(x,0) &=&g(x),\text{ }a\leq x\leq b  \notag
\end{eqnarray}%
will be decided in the later sections according to test problem. $U(x,t)$
and $V(x,t)$ are the unknown to be determined, $U_{t}\ $and $V_{t}$ are the
unsteady terms, $UU_{x}$ and $VV_{x}$ are the nonlinear terms and $U_{xx}$
and $V_{xx}$ are diffusive terms.

The approximate-analytical and numerical solutions of the CBE are in need
because precise analytical solution of the CBE does not exist in the
literature for wide class of boundary and initial conditions. The engineers
and scientists are studying the CBE to discover more properties of the CBE
due to its applicability in relevant fields. The approximate-analytical
methods such as decomposition method\cite{2001}, variational iteration method%
\cite{2005}, Tanh-method\cite{2006}, Adomian-pade technique \cite{2007},
differential transform method\cite{2011}, homotopy perturbation method\cite%
{2008a,2009a,2012b} are applied to find the approximate functional solutions
of the CBE.

Researcher have set up numerical methods to reveal more properties of the
CBE for specific solutions which does not found with analytical and
approximate analytical solutions. Episov presented numerical solution of the
one dimensional CBE and gave comparison of the results with the experimental
data\cite{1995,1995a}. Rashid et al. constructed the Fourier \
pseudo-spectral method for providing an approximate solutions for the CBE of
a set of initial and periodic boundary condition. A fully-implicit finite
difference method \cite{2013}, a composite numerical scheme based on finite
difference \cite{2014} and an implicit logarithmic finite difference method
have been presented for the CBE\cite{2014iki}. The Differential quadrature 
\cite{2012} \ and the generalized differential quadrature method have been
set up to obtain numerical solutions of the CBE. The collocation algorithms
based on chebyshev functions\cite{2008} and B-splines\cite{2011iki} are
written to study solutions of the CBE. The numerical solutions of CBE are
given by local discontinuous Galerkin method\cite{2011a}, the Galerkin $%
H^{2} $-Galerkin mixed method\cite{2012a} and Galerkin quadratic B-spline
finite element method \cite{2012iki}. Taylor Collocation-Extended Cubic
B-spline method is presented in the doctoral dissertation of Aksoy\cite%
{Aksoy}. The Modified cubic B-spline collocation method is presented by R.
C. Mittal and A. Tripathi\cite{2014uc}.

The B-splines are popular in science and engineering for solving the
differential equations. Main advantages of using the B-splines in the
numerical method is to provide simple and economical algorithms. The
trigonometric B-splines(TB) is an alternative base functions to the
well-known polynomial B-spline base functions. TB, introduced by I. J.
Schoenberg in 1964, is known as a non-polynomial B-splines including sine
function. These spline functions have started to use as an approximate
function for constructing the numerical methods to get numerical solutions
of the differential equations. The numerical methods, based on quadratic and
TB, for solving types of ordinary differential equations ware given in the
studies\cite{g1,g2,we,w1}. Very recently a collocation finite difference
scheme based on new TB is developed for the numerical solution of a
one-dimensional hyperbolic equation (wave equation) with non-local
conservation condition\cite{aa}. \ A new two-time level implicit technique
based on TB is proposed for the approximate solution of the nonclassical
diffusion problem with nonlocal boundary condition in the study\cite{ma}.

Different kinds of numerical methods have been developed to deal with
finding the solutions of the Burgers equation. \ It is known that solutions
of the Burgers type equations include the steep front and sharpness.
Therefore modelling of the solutions of the Burgers equation having sharp
behavior \ is of interest for the numerical analysist. \ In this paper, TB
are used to establish a collocation method for obtaining the numerical
solutions of the CBE. \ The fully-integration of the CBE is going to give \
an nonlinear algebraic equation. The time and space integrations have been
managed by employing the Crack-Nicolson method and the trigonometric
B-spline approximation respectively. After getting solution of the algebraic
system,\ the numerical solutions of the GBE will be found in terms of
combination of the cubic TB and the solutions of the algebraic equation. We
will present easy and simple algorithm with the use of the cubic TB in the
collocation method for solving CBE.

\section{Trigonometric Cubic B-spline Collocation Method}

\qquad Consider a uniform partition of the problem domain $[a=x_{0},b=x_{N}]$
at the knots $x_{i},i=0,...,N$ with mesh spacing $h=(b-a)/N.$ On this
partition together with additional knots $%
x_{-1},x_{0},x_{N+1},x_{N+2},x_{N+3}$ outside the problem domain, $CTB_{i}${}%
$(x)$ can be defined as

\begin{equation}
CTB_{i}(x)=\frac{1}{\theta }\left \{ 
\begin{tabular}{ll}
$\omega ^{3}(x_{i-2}),$ & $x\in \left[ x_{i-2},x_{i-1}\right] $ \\ 
$\omega (x_{i-2})(\omega (x_{i-2})\phi (x_{i})+\phi (x_{i+1})\omega
(x_{i-1}))+\phi (x_{i+2})\omega ^{2}(x_{i-1}),$ & $x\in \left[ x_{i-1},x_{i}%
\right] $ \\ 
$\omega (x_{i-2})\phi ^{2}(x_{i+1})+\phi (x_{i+2})(\omega (x_{i-1})\phi
(x_{i+1})+\phi (x_{i+2})\omega (x_{i})),$ & $x\in \left[ x_{i},x_{i+1}\right]
$ \\ 
$\phi ^{3}(x_{i+2}),$ & $x\in \left[ x_{i+1},x_{i+2}\right] $ \\ 
$0,$ & $\text{otherwise}$%
\end{tabular}%
\right.  \label{CB4}
\end{equation}%
where%
\begin{equation*}
\omega (x_{i})=\sin (\frac{x-x_{i}}{2}),\text{ }\phi (x_{i})=\sin (\frac{%
x_{i}-x}{2}),\text{ }\theta =\sin (\frac{h}{2})\sin (h)\sin (\frac{3h}{2}).
\end{equation*}

$CTB_{i}(x)$ are twice continuously differentiable piecewise trigonometric
cubic B-spline on the interval $[a,b]$. \ The iterative formula \bigskip 
\begin{equation}
T_{i}^{k}(x)=\frac{\sin (\frac{x-x_{i}}{2})}{\sin (\frac{x_{i+k-1}-x_{i}}{2})%
}T_{i}^{k-1}(x)+\frac{\sin (\frac{x_{i+k}-x}{2})}{\sin (\frac{x_{i+k}-x_{i+1}%
}{2})}T_{i+1}^{k-1}(x),\text{ }k=2,3,4,...  \label{CB5}
\end{equation}%
gives the cubic B-spline trigonometric functions starting with the
CTB-splines of order $1$%
\begin{equation*}
T_{i}^{1}(x)=\left \{ 
\begin{tabular}{c}
$1,\ x\in \lbrack x_{i},x_{i+1})$ \\ 
$0$ \  \ ,otherwise.%
\end{tabular}%
\right.
\end{equation*}

Each $CTB_{i}(x)$ is twice continuously differentiable and the values of $%
CTB_{i}(x),CTB_{i}^{^{\prime }}(x)$ and $CTB_{i}^{^{\prime \prime }}(x)$ at
the knots $x_{i}$ 's can be computed from Eq.(\ref{CB4}) as

\begin{equation*}
\begin{tabular}{l}
Table 1: Values of $B_{i}(x)$ and its principle two \\ 
derivatives at the knot points \\ 
\begin{tabular}{|l|l|l|l|}
\hline
& $T_{i}(x_{k})$ & $T_{i}^{\prime }(x_{k})$ & $T_{i}^{\prime \prime }(x_{k})$
\\ \hline
$x_{i-2}$ & $0$ & $0$ & $0$ \\ \hline
$x_{i-1}$ & $\sin ^{2}(\frac{h}{2})\csc \left( h\right) \csc (\frac{3h}{2})$
& $\frac{3}{4}\csc (\frac{3h}{2})$ & $\frac{3(1+3\cos (h))\csc ^{2}(\frac{h}{%
2})}{16\left[ 2\cos (\frac{h}{2})+\cos (\frac{3h}{2})\right] }$ \\ \hline
$x_{i}$ & $\frac{2}{1+2\cos (h)}$ & $0$ & $\frac{-3\cot ^{2}(\frac{3h}{2})}{%
2+4\cos (h)}$ \\ \hline
$x_{i+1}$ & $\sin ^{2}(\frac{h}{2})\csc \left( h\right) \csc (\frac{3h}{2})$
& $-\frac{3}{4}\csc (\frac{3h}{2})$ & $\frac{3(1+3\cos (h))\csc ^{2}(\frac{h%
}{2})}{16\left[ 2\cos (\frac{h}{2})+\cos (\frac{3h}{2})\right] }$ \\ \hline
$x_{i+2}$ & $0$ & $0$ & $0$ \\ \hline
\end{tabular}%
\end{tabular}%
\end{equation*}

$CTB_{i}(x)$ , $i=-1,...,N+1$ are a basis for the trigonometric spline
space. An approximate solution $U_{N}$ to the unknown $U$ is written in
terms of the expansion of the CTB as

\begin{equation}
U_{N}(x,t)=\sum_{i=-1}^{N+1}\delta _{i}CTB_{i}(x),\text{ }%
V_{N}(x,t)=\sum_{i=-1}^{N+1}\phi _{i}CTB_{i}(x)  \label{CBu1}
\end{equation}%
where $\delta _{i}$ are time dependent parameters to be determined from the
collocation points $x_{i},i=0,...,N$ and the boundary and initial
conditions. The first and second derivatives also can be defined by 
\begin{equation}
\begin{array}{c}
U_{N}^{\prime }(x,t)=\sum_{i=-1}^{N+1}\delta _{i}CTB_{i}^{\prime }(x),\text{ 
}V_{N}^{\prime }(x,t)=\sum_{i=-1}^{N+1}\phi _{i}CTB_{i}(x) \\ 
U_{N}^{\prime \prime }(x,t)=\sum_{i=-1}^{N+1}\delta _{i}CTB_{i}^{\prime
\prime }(x),\text{ }V_{N}^{\prime \prime }(x,t)=\sum_{i=-1}^{N+1}\phi
_{i}CTB_{i}(x)%
\end{array}
\label{CBu2}
\end{equation}

Using the Eq. (\ref{CBu1}), (\ref{CBu2}) and Table 1, we see that the nodal
values $U_{i},$ $V_{i}$, their first derivatives $U_{i}^{\prime },$ $%
V_{i}^{\prime }$ and second derivatives $U_{i}^{\prime \prime },$ $%
V_{i}^{\prime \prime }$ at the knots are given in terms of parameters by the
following relations%
\begin{equation}
\begin{array}{c}
\begin{tabular}{l}
$U_{i}=\alpha _{1}\delta _{i-1}+\alpha _{2}\delta _{i}+\alpha _{1}\delta
_{i+1}$ \\ 
$U_{i}^{\prime }=\beta _{1}\delta _{i-1}+\beta _{2}\delta _{i+1}$ \\ 
$U_{i}^{\prime \prime }=\gamma _{1}\delta _{i-1}+\gamma _{2}\delta
_{i}+\gamma _{1}\delta _{i+1}$%
\end{tabular}
\\ 
\begin{tabular}{l}
$V_{i}=\alpha _{1}\phi _{i-1}+\alpha _{2}\phi _{i}+\alpha _{1}\phi _{i+1}$
\\ 
$V_{i}^{\prime }=\beta _{1}\phi _{i-1}+\beta _{2}\phi _{i+1}$ \\ 
$V_{i}^{\prime \prime }=\gamma _{1}\phi _{i-1}+\gamma _{2}\phi _{i}+\gamma
_{1}\phi _{i+1}$%
\end{tabular}%
\end{array}
\label{CBu3}
\end{equation}%
where%
\begin{equation*}
\begin{array}{ll}
\alpha _{1}=\sin ^{2}(\frac{h}{2})\csc (h)\csc (\frac{3h}{2}) & \alpha _{2}=%
\dfrac{2}{1+2\cos (h)} \\ 
\beta _{1}=-\frac{3}{4}\csc (\frac{3h}{2}) & \beta _{2}=\frac{3}{4}\csc (%
\frac{3h}{2}) \\ 
\gamma _{1}=\dfrac{3((1+3\cos (h))\csc ^{2}(\frac{h}{2}))}{16(2\cos (\frac{h%
}{2})+\cos (\frac{3h}{2}))} & \gamma _{2}=-\dfrac{3\cot ^{2}(\frac{h}{2})}{%
2+4\cos (h)}%
\end{array}%
\end{equation*}

The Crank-Nicolson\ scheme is used to discretize time variables of the
unknown $U$ and $V$ in the Coupled Burger equation . we obtain the time
discretized form of the equation as%
\begin{equation}
\begin{array}{r}
\dfrac{U^{n+1}-U^{n}}{\Delta t}-\dfrac{U_{xx}^{n+1}+U_{xx}{}^{n}}{2}+k_{1}%
\dfrac{(UU_{x})^{n+1}+(UU_{x}^{n})^{n}}{2}+k_{2}\dfrac{%
(UV)_{x}^{n+1}+(UV)_{x}^{n}}{2}=0 \\ 
\dfrac{V^{n+1}-V^{n}}{\Delta t}-\dfrac{V_{xx}^{n+1}+V_{xx}{}^{n}}{2}+k_{1}%
\dfrac{(VV_{x})^{n+1}+(VV_{x}^{n})^{n}}{2}+k_{3}\dfrac{%
(UV)_{x}^{n+1}+(UV)_{x}^{n}}{2}=0%
\end{array}
\label{CB7}
\end{equation}%
where $U^{n+1}=U(x,t_{n}+\Delta t)$ and $V^{n+1}=V(x,t_{n}+\Delta t).$ The
nonlinear term $(UU_{x})^{n+1},$ $(VV_{x})^{n+1}$ and $(UV)_{x}^{n+1}$in Eq.
(\ref{CB7}) is linearized by using the following forms \cite{rubin}:

\begin{equation}
\begin{array}{ll}
(UU_{x})^{n+1} & =U^{n+1}U_{x}^{n}+U^{n}U_{x}^{n+1}-U^{n}U_{x}^{n} \\ 
(VV_{x})^{n+1} & =V^{n+1}V_{x}^{n}+V^{n}V_{x}^{n+1}-V^{n}V_{x}^{n} \\ 
(UV)_{x}^{n+1} & =(U_{x}V)^{n+1}+(UV_{x})^{n+1} \\ 
& 
=U_{x}^{n+1}V^{n}+U_{x}^{n}V^{n+1}-U_{x}^{n}V^{n}+U^{n+1}V_{x}^{n}+U^{n}V_{x}^{n+1}-U^{n}V_{x}^{n}%
\end{array}
\label{Rubin}
\end{equation}

Substitution the expressions (\ref{CBu3}) into (\ref{CB7}) and evaluating
the resulting \ equations at the knots yields the system of the
fully-discretized equations

\begin{equation}
\nu _{m1}\delta _{m-1}^{n+1}+\nu _{m2}\phi _{m-1}^{n+1}+\nu _{m3}\delta
_{m}^{n+1}+\nu _{m4}\phi _{m}^{n+1}+\nu _{m5}\delta _{m+1}^{n+1}+\nu
_{m6}\phi _{m+1}^{n+1}=\nu _{m7}\delta _{m-1}^{n}+\nu _{m8}\delta
_{m}^{n}+\nu _{m9}\delta _{m+1}^{n}  \label{sis1}
\end{equation}

and%
\begin{equation}
\nu _{m10}\delta _{m-1}^{n+1}+\nu _{m11}\phi _{m-1}^{n+1}+\nu _{m12}\delta
_{m}^{n+1}+\nu _{m13}\phi _{m}^{n+1}+\nu _{m14}\delta _{m+1}^{n+1}+\nu
_{m15}\phi _{m+1}^{n+1}=\nu _{m7}\phi _{m-1}^{n}+\nu _{m8}\phi _{m}^{n}+\nu
_{m9}\phi _{m+1}^{n}  \label{sis2}
\end{equation}

where%
\begin{equation*}
\begin{array}{l}
\nu _{m1}=\left( \dfrac{2}{\Delta t}+k_{1}K_{2}+k_{2}L_{2}\right) \alpha
_{1}+\left( k_{1}K_{1}+k_{2}L_{1}\right) \beta _{1}-\gamma _{1} \\ 
\nu _{m2}=\left( k_{1}K_{2}\right) \alpha _{1}+\left( k_{2}K_{1}\right)
\beta _{1} \\ 
\nu _{m3}=\left( \dfrac{2}{\Delta t}+k_{1}K_{2}+k_{2}L_{2}\right) \alpha
_{2}-\gamma _{2} \\ 
\nu _{m4}=\left( k_{1}K_{2}\right) \alpha _{2} \\ 
\nu _{m5}=\left( \dfrac{2}{\Delta t}+k_{1}K_{2}+k_{2}L_{2}\right) \alpha
_{1}-\left( k_{1}K_{1}+k_{2}L_{1}\right) \beta _{2}-\gamma _{1} \\ 
\nu _{m6}=\left( k_{1}K_{2}\right) \alpha _{1}-\left( k_{2}K_{1}\right)
\beta _{2} \\ 
\\ 
\nu _{m7}=\frac{2}{\Delta t}\alpha _{1}+\gamma _{1} \\ 
\nu _{m8}=\frac{2}{\Delta t}\alpha _{2}+\gamma _{2} \\ 
\nu _{m9}=\frac{2}{\Delta t}\alpha _{1}+\gamma _{1} \\ 
\\ 
\nu _{m10}=\left( k_{3}L_{2}\right) \alpha _{1}+\left( k_{3}L_{1}\right)
\beta _{1} \\ 
\nu _{m11}=\left( \frac{2}{\Delta t}+k_{1}L_{2}+k_{3}K_{2}\right) \alpha
_{1}+\left( k_{1}L_{1}+k_{3}K_{1}\right) \beta _{1}-\gamma _{1} \\ 
\nu _{m12}=\left( k_{3}L_{2}\right) \alpha _{2} \\ 
\nu _{m12}=\left( \frac{2}{\Delta t}+k_{1}L_{2}+k_{3}K_{2}\right) \alpha
_{2}-\gamma _{2} \\ 
\nu _{m14}=\left( k_{3}L_{2}\right) \alpha _{1}-\left( k_{3}L_{1}\right)
\beta _{2} \\ 
\nu _{m15}=\left( \frac{2}{\Delta t}+k_{1}L_{2}+k_{3}K_{2}\right) \alpha
_{1}-\left( k_{1}L_{1}+k_{3}K_{1}\right) \beta _{2}-\gamma _{1}%
\end{array}%
\end{equation*}%
\begin{eqnarray*}
&&%
\begin{array}{cc}
K_{1}=\alpha _{1}\delta _{m-1}^{n}+\alpha _{2}\delta _{m}^{n}+\alpha
_{3}\delta _{m+1}^{n} & L_{1}=\alpha _{1}\phi _{m-1}^{n}+\alpha _{2}\phi
_{m}^{n}+\alpha _{3}\phi _{m+1}^{n}%
\end{array}
\\
&&%
\begin{array}{cc}
K_{2}=\beta _{1}\delta _{m-1}^{n}+\beta _{2}\delta _{m+1}^{n} & L_{2}=\beta
_{1}\phi _{m-1}^{n}+\beta _{2}\phi _{m+1}^{n}%
\end{array}%
\end{eqnarray*}%
The system with (\ref{sis1}) and (\ref{sis2}) can be expressed in the
following matrices system;%
\begin{equation}
\mathbf{Ad}^{n+1}=\mathbf{Bd}^{n}  \label{CB10}
\end{equation}%
where%
\begin{equation*}
\mathbf{A=}%
\begin{bmatrix}
\nu _{m1} & \nu _{m2} & \nu _{m3} & \nu _{m4} & \nu _{m5} & \nu _{m6} &  & 
&  &  \\ 
\nu _{m10} & \nu _{m11} & \nu _{m12} & \nu _{m13} & \nu _{m14} & \nu _{m15}
&  &  &  &  \\ 
&  & \nu _{m1} & \nu _{m2} & \nu _{m3} & \nu _{m4} & \nu _{m5} & \nu _{m6} & 
&  \\ 
&  & \nu _{m10} & \nu _{m11} & \nu _{m12} & \nu _{m13} & \nu _{m14} & \nu
_{m15} &  &  \\ 
&  &  & \ddots & \ddots & \ddots & \ddots & \ddots & \ddots &  \\ 
&  &  &  & \nu _{m1} & \nu _{m2} & \nu _{m3} & \nu _{m4} & \nu _{m5} & \nu
_{m6} \\ 
&  &  &  & \nu _{m10} & \nu _{m11} & \nu _{m12} & \nu _{m13} & \nu _{m14} & 
\nu _{m15}%
\end{bmatrix}%
\end{equation*}%
\begin{equation*}
\mathbf{B=}%
\begin{bmatrix}
\nu _{m7} & 0 & \nu _{m8} & 0 & \nu _{m9} & 0 &  &  &  &  \\ 
0 & \nu _{m7} & 0 & \nu _{m8} & 0 & \nu _{m9} &  &  &  &  \\ 
&  & \nu _{m7} & 0 & \nu _{m8} & 0 & \nu _{m9} & 0 &  &  \\ 
&  & 0 & \nu _{m7} & 0 & \nu _{m8} & 0 & \nu _{m9} &  &  \\ 
&  &  & \ddots & \ddots & \ddots & \ddots & \ddots & \ddots &  \\ 
&  &  &  & \nu _{m7} & 0 & \nu _{m8} & 0 & \nu _{m9} & 0 \\ 
&  &  &  & 0 & \nu _{m7} & 0 & \nu _{m8} & 0 & \nu _{m9}%
\end{bmatrix}%
\end{equation*}

\qquad The system (\ref{CB10}) consist of $2N+2$ linear equation in $2N+6$
unknown parameters $\mathbf{d}^{n+1}=(\delta _{-1}^{n+1},\phi
_{-1}^{n+1},\delta _{0}^{n+1},\phi _{0}^{n+1},\ldots ,\delta
_{n+1}^{n+1},\phi _{n+1}^{n+1},)$. To obtain a unique solution an additional
four constraints\ are needed. By imposing the Dirichlet boundary conditions
this will lead us to the following relations;%
\begin{equation}
\begin{array}{r}
\delta _{-1}=(f_{1}(a,t)-\alpha _{2}\delta _{0}-\alpha _{1}\delta
_{1})/\alpha _{1} \\ 
\phi _{-1}=(g_{1}(a,t)-\alpha _{2}\phi _{0}-\alpha _{1}\phi _{1})/\alpha _{1}
\\ 
\delta _{N+1}=(f_{2}(b,t)-\alpha _{1}\delta _{N-1}-\alpha _{2}\delta
_{N})/\alpha _{1} \\ 
\phi _{N+1}=(g_{2}(b,t)-\alpha _{1}\phi _{N-1}+\alpha _{2}\phi _{N})/\alpha
_{1}%
\end{array}
\label{CB11}
\end{equation}

Elimination of the parameters $\delta _{-1},$ $\phi _{-1},$ $\delta _{N+1},$ 
$\phi _{N+1}$ from the system (\ref{CB10}) gives us a solvable system. Once
the initial parameters are determined, parameters $\delta
_{i}^{n},i=-1...N+1 $ can be found from the iterative system (\ref{CB10}),
so that using the equations solutions of the equations at the defined time
steps can be found iteratively.

\section{The Initial State}

\qquad Initial parameters $\delta _{-1}^{0},$ $\phi _{-1}^{0},$ $\delta
_{0}^{0},$ $\phi _{0}^{0},\ldots ,\delta _{N+1}^{0},$ $\phi _{N+1}^{0}$ can
be determined from the initial condition and first space derivative of the
initial conditions at the boundaries as the following:%
\begin{equation}
\begin{array}{l}
U^{0}(a,0)=\alpha _{1}\delta _{-1}^{0}+\alpha _{2}\delta _{0}^{0}+\alpha
_{1}\delta _{1}^{0} \\ 
U^{0}(x_{m},0)=\alpha _{1}\delta _{m-1}^{0}+\alpha _{2}\delta
_{m}^{0}+\alpha _{1}\delta _{m+1}^{0},\text{ \  \ }m=1,2,...,N-1 \\ 
U^{0}(b,0)=U_{N}^{0}=\alpha _{1}\delta _{N-1}^{0}+\alpha _{2}\delta
_{N}^{0}+\alpha _{1}\delta _{N+1}^{0}%
\end{array}
\label{BU1}
\end{equation}%
and%
\begin{equation}
\begin{array}{l}
V^{0}(a,0)=\alpha _{1}\phi _{-1}^{0}+\alpha _{2}\phi _{0}^{0}+\alpha
_{1}\phi _{1}^{0} \\ 
V^{0}(x_{m},0)=\alpha _{1}\phi _{m-1}^{0}+\alpha _{2}\phi _{m}^{0}+\alpha
_{1}\phi _{m+1}^{0},\text{ \  \ }m=1,2,...,N-1 \\ 
V^{0}(b,0)=\alpha _{1}\phi _{N-1}^{0}+\alpha _{2}\phi _{N}^{0}+\alpha
_{1}\phi _{N+1}^{0}%
\end{array}
\label{BV1}
\end{equation}

The system (\ref{BU1}) consists $N+1$ equations and $N+3$ unknown, so we
have to eliminate $\delta _{-1}^{0}$ and $\delta _{N+1}^{0}$ for solving
this system using following derivatives conditions%
\begin{equation*}
\delta _{-1}^{0}=\left( U_{0}^{\prime }-\beta _{2}\delta _{1}^{0}\right)
/\beta _{1},\text{ }\delta _{N+1}^{0}=\left( U_{N}^{\prime }-\beta
_{1}\delta _{N-1}^{0}\right) /\beta _{2},
\end{equation*}%
and if the equations system is rearranged for the above conditions, then
following form is obtained%
\begin{equation*}
\left[ 
\begin{tabular}{ccccccc}
$\smallskip \alpha _{2}$ & $2\alpha _{1}$ &  &  &  &  &  \\ 
$\alpha _{1}$ & $\alpha _{2}$ & $\alpha _{1}$ &  &  &  &  \\ 
&  &  & $\ddots $ &  &  &  \\ 
&  &  &  & $\alpha _{1}\smallskip $ & $\alpha _{2}$ & $\alpha _{1}$ \\ 
&  &  &  &  & $2\alpha _{1}$ & $\alpha _{2}$%
\end{tabular}%
\right] \left[ 
\begin{tabular}{c}
$\delta _{0}^{0}\smallskip $ \\ 
$\delta _{1}^{0}$ \\ 
$\vdots $ \\ 
$\delta _{N-1}^{0}\smallskip $ \\ 
$\delta _{N}^{0}$%
\end{tabular}%
\right] =\left[ 
\begin{tabular}{c}
$U_{0}^{\prime }-\tfrac{\alpha _{1}}{\beta _{1}}U_{0}^{\prime }\smallskip $
\\ 
$U_{1}^{\prime }$ \\ 
$\vdots $ \\ 
$U_{N-1}^{\prime }\smallskip $ \\ 
$U_{N}^{\prime }-\tfrac{\alpha _{1}}{\beta _{1}}U_{N}^{\prime }$%
\end{tabular}%
\right]
\end{equation*}%
which can also be solved using a variant of the Thomas algorithm. With the
same way, from the system (\ref{BV1}), $\phi _{-1}^{0}$ and $\phi _{N+1}^{0}$
can be eliminated using%
\begin{equation*}
\phi _{-1}^{0}=\left( V_{0}^{\prime }-\beta _{2}\phi _{1}^{0}\right) /\beta
_{1},\text{ }\phi _{N+1}^{0}=\left( V_{N}^{\prime }-\beta _{1}\phi
_{N-1}^{0}\right) /\beta _{2},
\end{equation*}%
conditions and the following three bounded matrix is obtained.%
\begin{equation*}
\left[ 
\begin{tabular}{ccccccc}
$\smallskip \alpha _{2}$ & $2\alpha _{1}$ &  &  &  &  &  \\ 
$\alpha _{1}$ & $\alpha _{2}$ & $\alpha _{1}$ &  &  &  &  \\ 
&  &  & $\ddots $ &  &  &  \\ 
&  &  &  & $\alpha _{1}\smallskip $ & $\alpha _{2}$ & $\alpha _{1}$ \\ 
&  &  &  &  & $2\alpha _{1}$ & $\alpha _{2}$%
\end{tabular}%
\right] \left[ 
\begin{tabular}{c}
$\phi _{0}^{0}\smallskip $ \\ 
$\phi _{1}^{0}$ \\ 
$\vdots $ \\ 
$\phi _{N-1}^{0}\smallskip $ \\ 
$\phi _{N}^{0}$%
\end{tabular}%
\right] =\left[ 
\begin{tabular}{c}
$V_{0}^{\prime }-\tfrac{\alpha _{1}}{\beta _{1}}V_{0}^{\prime }\smallskip $
\\ 
$V_{1}^{\prime }$ \\ 
$\vdots $ \\ 
$V_{N-1}^{\prime }\smallskip $ \\ 
$V_{N}^{\prime }-\tfrac{\alpha _{1}}{\beta _{1}}V_{N}^{\prime }$%
\end{tabular}%
\right]
\end{equation*}

\section{Numerical Tests}

\qquad In this section, numerical results of the three test problems will be
presented for the coupled Burgers equation. The accuracy of suggested method
problem will be shown by calculating the error norm%
\begin{equation*}
L_{\infty }=\left \vert U-U_{N}\right \vert _{\infty }=\max \limits_{j}\left
\vert U_{j}-(U_{N})_{j}^{n}\right \vert
\end{equation*}
The obtained results will compare with results of \cite{2011iki}, \cite{2012}%
, \cite{2014uc} and \cite{Aksoy}.

\textbf{Problem 1) }Consider the coupled Burgers equation (\ref{CB1}) with
the following initial and boundary conditions%
\begin{equation*}
U(x,0)=\sin (x),\text{ }V(x,0)=\sin (x)
\end{equation*}%
and%
\begin{equation*}
U(-\pi ,t)=U(\pi ,t)=V(-\pi ,t)=V(\pi ,t)=0
\end{equation*}%
The exact solution is%
\begin{equation*}
U(x,t)=V(x,t)=e^{-t}\sin (x)
\end{equation*}

We compute the numerical solutions using the selected values $k_{1}=-2,$ $%
k_{2}=1$ and $k_{3}=1$ with different values of time step length $\Delta t.$%
In our first computation, we take $t=0.1,$ $\Delta t=0.001$ while the number
of partition $N$ changes$.$ The corresponding results are presented in Table
2 a. In our computation, we compute the maximum absolute errors at time
level $t=1$ for the parameters with different decreasing values of $t$. The
corresponding results are reported in Table 2 b. In both computations, the
results are same for $U(x,t)$ and $V(x,t)$ because of symmetric initial and
boundary conditions. And also we correspond the obtained numerical solutions
by different settings of parameters, specifically for those taken by \cite%
{2014uc} in Table 2 c for $N=50$, $\Delta t=0.01$ and increasing $t.$ And
also in Table 2, we present the space rate of convergence for $t=3$ which is
clearly of second order.

\begin{equation*}
\begin{tabular}{|l|}
\hline
Table 2 a: $L_{\infty }$ Error norms for $t=0.1,$ $\Delta t=0.001,$ $%
U(x,t)=V(x,t)$ \\ \hline
\begin{tabular}{ccccc}
& Present & \cite{2011iki} & \cite{Aksoy}, ($\lambda =0)$ & \cite{Aksoy}
(Various $\lambda $) \\ 
$N=200$ & $0.69699\times 10^{-5}$ & $0.745\times 10^{-5}$ & $0.74326\times
10^{-5}$ & $0.00079\times 10^{-5}(\lambda =-1.640\times 10^{-4})$ \\ 
$N=400$ & $0.17367\times 10^{-5}$ & $0.186\times 10^{-5}$ & $0.18534\times
10^{-5}$ & $0.00006\times 10^{-5}(\lambda =-4.087\times 10^{-5})$%
\end{tabular}
\\ \hline
\end{tabular}%
\end{equation*}%
\begin{equation*}
\begin{tabular}{|l|}
\hline
Table 2 b: $L_{\infty }$ Error norms for $t=1,$ $N=400,$ $U(x,t)=V(x,t)$ \\ 
\hline
\begin{tabular}{cccccc}
& Present & \cite{Aksoy}, ($\lambda =0)$ & \cite{Aksoy} (Various $\lambda $)
& \cite{2011iki} & \cite{Rashid} \\ 
$\Delta t=0.01$ & $0.40261\times 10^{-5}$ & $1.08691\times 10^{-5}$ & $%
0.00131\times 10^{-5}(\lambda =-5.896\times 10^{-5})$ &  &  \\ 
$\Delta t=0.001$ & $0.70610\times 10^{-5}$ & $1.10393\times 10^{-5}$ & $%
0.00036\times 10^{-5}(\lambda =-5.992\times 10^{-5})$ & $0.756\times 10^{-5}$
& $0.116\times 10^{-4}$%
\end{tabular}
\\ \hline
\end{tabular}%
\end{equation*}%
\begin{equation*}
\begin{tabular}{|l|}
\hline
Table 2 c: $L_{\infty }$ Error norms for $\Delta t=0.01,$ $N=50,$ $%
U(x,t)=V(x,t)$ different $t.$ \\ \hline
\begin{tabular}{ccccc}
& Present & \cite{2012} & \cite{2012iki} & \cite{2014uc} \\ 
$t=0.5$ & $3.7144\times 10^{-4}$ & $1.51688\times 10^{-4}$ & $2.26627\times
10^{-5}$ & $1.103080984\times 10^{-4}$ \\ 
$t=1.0$ & $4.5072\times 10^{-4}$ & $1.83970\times 10^{-4}$ & $1.46179\times
10^{-5}$ & $1.336880384\times 10^{-4}$ \\ 
$t=2.0$ & $3.3183\times 10^{-4}$ & $1.35250\times 10^{-4}$ & $0.73805\times
10^{-5}$ & $9.818252567\times 10^{-5}$ \\ 
$t=3.0$ & $1.8322\times 10^{-4}$ & $7.46014\times 10^{-5}$ & $0.40272\times
10^{-5}$ & $1.029870405\times 10^{-5}$%
\end{tabular}
\\ \hline
\end{tabular}%
\end{equation*}%
\begin{equation*}
\begin{tabular}{|l|}
\hline
Table 2 d: The rate of convergence for $\Delta t=0.01$ and $\Delta t=0.0001$
respectively \\ \hline
\begin{tabular}{ccc}
$\Delta t=0.01$ & Present & order \\ 
$N=50$ & $1.8322\times 10^{-4}$ &  \\ 
$N=100$ & $4.4857\times 10^{-5}$ & $2.0302$ \\ 
$N=150$ & $1.9232\times 10^{-5}$ & $2.0887$ \\ 
$N=200$ & $1.0274\times 10^{-5}$ & $2.1793$ \\ 
$N=250$ & $6.1264\times 10^{-6}$ & $2.3170$%
\end{tabular}%
\begin{tabular}{ccc}
$\Delta t=0.0001$ & Present ($p=1$) & order \\ 
$N=50$ & $1.8446\times 10^{-4}$ &  \\ 
$N=100$ & $4.6102\times 10^{-5}$ & $2.0004$ \\ 
$N=150$ & $2.0476\times 10^{-5}$ & $2.0016$ \\ 
$N=200$ & $1.1519\times 10^{-5}$ & $1.9998$ \\ 
$N=250$ & $7.3709\times 10^{-5}$ & $2.0007$%
\end{tabular}
\\ \hline
\end{tabular}%
\end{equation*}%
The corresponding graphical illustrations are presented in Figures 2 for $%
k_{1}=-2,$ $k_{2}=1$, $k_{3}=1,$ $N=400$ and $\Delta t=0.001$ at different $%
t $ for best parameter $p=0.0002166$. In Figure 2-4, computed solutions of $V
$ different time levels for $k_{1},$ $k_{2}$ fixed, $k_{1},$ $k_{3}$ and $%
k_{2},$ $k_{3}$ fixed respectively.%
\begin{equation*}
\begin{array}{c}
\FRAME{itbpF}{3.7075in}{3.7075in}{0in}{}{}{fig3.bmp}{\special{language
"Scientific Word";type "GRAPHIC";maintain-aspect-ratio TRUE;display
"USEDEF";valid_file "F";width 3.7075in;height 3.7075in;depth
0in;original-width 3.6599in;original-height 3.6599in;cropleft "0";croptop
"1";cropright "1";cropbottom "0";filename 'Fig3.bmp';file-properties
"XNPEU";}} \\ 
\text{Figure 2: Computed solutions of }U\text{ Problem 1 for }t=3,\text{ }%
k_{1}=-2,\text{ }k_{2}=1,\text{ }k_{3}=-8%
\end{array}%
\end{equation*}%
\begin{equation*}
\begin{array}{c}
\FRAME{itbpF}{3.7075in}{3.7075in}{0in}{}{}{fig3.bmp}{\special{language
"Scientific Word";type "GRAPHIC";maintain-aspect-ratio TRUE;display
"USEDEF";valid_file "F";width 3.7075in;height 3.7075in;depth
0in;original-width 3.6599in;original-height 3.6599in;cropleft "0";croptop
"1";cropright "1";cropbottom "0";filename 'Fig3.bmp';file-properties
"XNPEU";}} \\ 
\text{Figure 3: Computed solutions of }U\text{ Problem 1 for }t=3,\text{ }%
k_{1}=-2,\text{ }k_{2}=1,\text{ }k_{3}=-4%
\end{array}%
\end{equation*}%
\begin{equation*}
\begin{array}{c}
\FRAME{itbpF}{3.7075in}{3.7075in}{0in}{}{}{fig4.bmp}{\special{language
"Scientific Word";type "GRAPHIC";maintain-aspect-ratio TRUE;display
"USEDEF";valid_file "F";width 3.7075in;height 3.7075in;depth
0in;original-width 3.6599in;original-height 3.6599in;cropleft "0";croptop
"1";cropright "1";cropbottom "0";filename 'Fig4.bmp';file-properties
"XNPEU";}} \\ 
\text{Figure 4: Computed solutions of }U\text{ Problem 1 for }t=3,\text{ }%
k_{1}=-2,\text{ }k_{2}=1,\text{ }k_{3}=0%
\end{array}%
\end{equation*}%
\textbf{\ Problem 2) }Numerical solutions of considered coupled Burgers'
equations are obtained for $k_{1}=2$ with different values of $k_{2}$ and $%
k_{3}$ at different time levels. In this situation the exact solution is%
\begin{equation*}
\begin{array}{l}
U(x,t)=a_{0}-2A(\dfrac{2k_{2}-1}{4k_{2}k_{3}-1})\tanh (A(x-2At)) \\ 
V(x,t)=a_{0}(\dfrac{2k_{3}-1}{2k_{2}-1})-2A(\dfrac{2k_{2}-1}{4k_{2}k_{3}-1}%
)\tanh (A(x-2At))%
\end{array}%
\end{equation*}

Thus, the initial and boundary conditions are taken from the exact solution
is%
\begin{equation*}
\begin{array}{l}
U(x,0)=a_{0}-2A(\dfrac{2k_{2}-1}{4k_{2}k_{3}-1})\tanh (Ax) \\ 
V(x,0)=a_{0}(\dfrac{2k_{3}-1}{2k_{2}-1})-2A(\dfrac{2k_{2}-1}{4k_{2}k_{3}-1}%
)\tanh (Ax)%
\end{array}%
\end{equation*}

Thus, the initial and boundary conditions are extracted from the exact
solution. Where $a_{0}=0.05$ and $A=\dfrac{1}{2}(\dfrac{a_{0}(4k_{2}k_{3}-1)%
}{2k_{2}-1}).$ The numerical solutions have been computed for the domain $%
x\in \lbrack -10,10],$ $\Delta t=0.01$ and number of partition $N=100.$ The
maximum error norms have been computed and compared in Tables 3 a-3 b for $%
t=0.5$ and $t=1.0$ with those available in the literature \cite{2011iki}, 
\cite{2008} and \cite{Rashid}. For $N=21$ the maximum error norms have been
computed and compared in Tables 3 c-3 d for $t=0.5,$ $t=1.0$ and $t=3.0$
with those available in the literature \cite{2008}, \cite{Rashid}, \cite%
{2012} and \cite{2014uc}.%
\begin{equation*}
\begin{tabular}{|l|}
\hline
Table 3 a: Maximum error norms for $U(x,t),$ $k_{1}=2,$ $\Delta t=0.01,$ $%
N=100$ \\ \hline
\begin{tabular}{ccccccc}
$t$ & $k_{2}$ & $k_{3}$ & Present & \cite{2011iki} & \cite{2008} & \cite%
{Rashid} \\ 
$0.5$ & $0.1$ & $0.3$ & $0.4707\times 10^{-4}$ & $0.4167\times 10^{-4}$ & $%
0.438\times 10^{-4}$ & $0.9619\times 10^{-3}$ \\ 
& $0.3$ & $0.03$ & $0.2709\times 10^{-4}$ & $0.4590\times 10^{-4}$ & $%
0.458\times 10^{-4}$ & $0.4310\times 10^{-3}$ \\ 
$1.0$ & $0.1$ & $0.3$ & $0.2831\times 10^{-4}$ & $0.8258\times 10^{-4}$ & $%
0.866\times 10^{-4}$ & $0.1152\times 10^{-2}$ \\ 
& $0.3$ & $0.03$ & $0.4988\times 10^{-4}$ & $0.9182\times 10^{-4}$ & $%
0.916\times 10^{-4}$ & $0.1268\times 10^{-2}$%
\end{tabular}
\\ \hline
\end{tabular}%
\end{equation*}%
\begin{equation*}
\begin{tabular}{|l|}
\hline
Table 3 b: Maximum error norms for $V(x,t),$ $k_{1}=2,$ $\Delta t=0.01,$ $%
N=100$ \\ \hline
\begin{tabular}{ccccccc}
$t$ & $k_{2}$ & $k_{3}$ & Present & \cite{2011iki} & \cite{2008} & \cite%
{Rashid} \\ 
$0.5$ & $0.1$ & $0.3$ & $0.1247\times 10^{-4}$ & $0.1480\times 10^{-3}$ & $%
0.499\times 10^{-4}$ & $0.3332\times 10^{-3}$ \\ 
& $0.3$ & $0.03$ & $0.7641\times 10^{-4}$ & $0.5729\times 10^{-3}$ & $%
0.181\times 10^{-3}$ & $0.1148\times 10^{-2}$ \\ 
$1.0$ & $0.1$ & $0.3$ & $0.2474\times 10^{-4}$ & $0.4770\times 10^{-4}$ & $%
0.992\times 10^{-4}$ & $0.1162\times 10^{-2}$ \\ 
& $0.3$ & $0.03$ & $0.1523\times 10^{-4}$ & $0.3617\times 10^{-3}$ & $%
0.362\times 10^{-3}$ & $0.1638\times 10^{-2}$%
\end{tabular}
\\ \hline
\end{tabular}%
\end{equation*}%
\begin{equation*}
\begin{tabular}{|l|}
\hline
Table 3 c: Maximum error norms for $U(x,t)$ in Problem 2 $k_{1}=2,$ $\Delta
t=0.01,$ $N=21$ \\ \hline
\begin{tabular}{cccccccc}
$t$ & $k_{2}$ & $k_{3}$ & Present & \cite{2008} & \cite{Rashid} & \cite{2012}
& \cite{2014uc} \\ 
$0.5$ & $0.1$ & $0.3$ & $0.13145\times 10^{-2}$ & $0.144\times 10^{-2}$ & $%
0.9619\times 10^{-3}$ & $0.4173\times 10^{-4}$ & $0.4189217417\times 10^{-4}$
\\ 
& $0.3$ & $0.03$ & $0.15686\times 10^{-2}$ & $0.668\times 10^{-3}$ & $%
0.4310\times 10^{-3}$ & $0.4585\times 10^{-4}$ & $0.4584830094\times 10^{-4}$
\\ 
$1.0$ & $0.1$ & $0.3$ & $0.25126\times 10^{-2}$ & $0.127\times 10^{-2}$ & $%
0.1153\times 10^{-2}$ & $0.8275\times 10^{-4}$ & $0.8269641708\times 10^{-4}$
\\ 
& $0.3$ & $0.03$ & $0.29666\times 10^{-2}$ & $0.130\times 10^{-2}$ & $%
0.1268\times 10^{-2}$ & $0.9167\times 10^{-4}$ & $0.9147335667\times 10^{-4}$
\\ 
$3.0$ & $0.1$ & $0.3$ & $0.68877\times 10^{-2}$ &  &  & $0.2408\times
10^{-3} $ & $0.2401202768\times 10^{-3}$ \\ 
& $0.1$ & $0.03$ & $0.70013\times 10^{-2}$ &  &  & $0.2747\times 10^{-3}$ & $%
0.2704203611\times 10^{-3}$%
\end{tabular}
\\ \hline
\end{tabular}%
\end{equation*}%
\begin{equation*}
\begin{tabular}{|l|}
\hline
Table 3 d: Maximum error norms for $V(x,t)$ in Problem 2 $k_{1}=2,$ $\Delta
t=0.01,$ $N=21$ \\ \hline
\begin{tabular}{ccccccccc}
$t$ & $k_{2}$ & $k_{3}$ & Present &  & \cite{2008} & \cite{Rashid} & \cite%
{2012} & \cite{2014uc} \\ 
$0.5$ & $0.1$ & $0.3$ & $0.25322\times 10^{-5}$ & $0.78143\times 10^{-3}$ & $%
0.542\times 10^{-3}$ & $0.3332\times 10^{-3}$ & $0.5418\times 10^{-4}$ & $%
0.9094743099\times 10^{-5}$ \\ 
& $0.3$ & $0.03$ & $0.10429\times 10^{-4}$ & $0.29775\times 10^{-2}$ & $%
0.120\times 10^{-2}$ & $0.1148\times 10^{-2}$ & $0.2826\times 10^{-4}$ & $%
0.2482188810\times 10^{-4}$ \\ 
$1.0$ & $0.1$ & $0.3$ & $0.25602\times 10^{-5}$ & $0.14757\times 10^{-2}$ & $%
0.129\times 10^{-2}$ & $0.1162\times 10^{-2}$ & $0.1074\times 10^{-3}$ & $%
0.1696286567\times 10^{-4}$ \\ 
& $0.3$ & $0.03$ & $0.10512\times 10^{-4}$ & $0.57049\times 10^{-2}$ & $%
0.235\times 10^{-2}$ & $0.1638\times 10^{-2}$ & $0.5673\times 10^{-4}$ & $%
0.4965329678\times 10^{-4}$ \\ 
$3.0$ & $0.1$ & $0.3$ & $0.25980\times 10^{-5}$ & $0.39769\times 10^{-2}$ & 
&  & $0.3119\times 10^{-3}$ & $0.4505480184\times 10^{-4}$ \\ 
& $0.1$ & $0.03$ & $0.10542\times 10^{-4}$ & $0.80324\times 10^{-2}$ &  &  & 
$0.1663\times 10^{-3}$ & $0.1498311672\times 10^{-4}$%
\end{tabular}
\\ \hline
\end{tabular}%
\end{equation*}

\begin{equation*}
\begin{tabular}{c}
$%
\begin{array}{cc}
\FRAME{itbpF}{2.7224in}{2.3134in}{0in}{}{}{fig5a.bmp}{\special{language
"Scientific Word";type "GRAPHIC";maintain-aspect-ratio TRUE;display
"USEDEF";valid_file "F";width 2.7224in;height 2.3134in;depth
0in;original-width 2.6801in;original-height 2.2736in;cropleft "0";croptop
"1";cropright "1";cropbottom "0";filename 'Fig5a.bmp';file-properties
"XNPEU";}} & \FRAME{itbpF}{2.7155in}{2.2667in}{0in}{}{}{fig5b.bmp}{\special%
{language "Scientific Word";type "GRAPHIC";maintain-aspect-ratio
TRUE;display "USEDEF";valid_file "F";width 2.7155in;height 2.2667in;depth
0in;original-width 2.674in;original-height 2.2269in;cropleft "0";croptop
"1";cropright "1";cropbottom "0";filename 'Fig5b.bmp';file-properties
"XNPEU";}}%
\end{array}%
$ \\ 
Figure 5: Numerical Solutions for $U(x,t)$ and $V(x,t)$, $N=21,$ $\Delta
t=0.001,$ $t=1,$ $0\leq x\leq 1,$ $k_{2}=0.1$ and $k_{3}=0.3$%
\end{tabular}%
\end{equation*}%
\textbf{Problem 3)\ }Consider the Coupled Burger equation (\ref{CB1}) with
the following initial conditions

\begin{equation*}
U(x,0)=\left \{ 
\begin{array}{cc}
\sin (2\pi x), & x\in \lbrack 0,0.5] \\ 
0, & x\in (0.5,1]%
\end{array}%
\right.
\end{equation*}

\begin{equation*}
V(x,0)=\left \{ 
\begin{array}{cc}
0, & x\in \lbrack 0,0.5] \\ 
-\sin (2\pi x), & x\in (0.5,1]%
\end{array}%
\right.
\end{equation*}%
and zero boundary conditions. In the Problem 3, the solutions have been
carried out on $x\in \lbrack 0,1]$ with $\Delta t=0.001$ and number of
partitions as $50$. Maximum values of $u$ and $v$ at different time levels
for $k_{2}=k_{3}=10$ have been given in Table 4a and 4 b, while the Tables 4
c and 4 d represent the maximum values for $k_{2}=k_{3}=100.$%
\begin{equation*}
\begin{tabular}{|l|}
\hline
Table 4 a: Maximum values of $U$ at different time levels for $%
k_{2}=k_{3}=10 $ \\ \hline
\begin{tabular}{ccccc}
$t$ & Present & \cite{2011iki} & \cite{2014uc} & at point \\ 
$0.1$ & $0.142427$ & $0.14456$ & $0.144491495800$ & $0.58$ \\ 
$0.2$ & $0.051716$ & $0.05237$ & $0.052356151890$ & $0.54$ \\ 
$0.3$ & $0.019087$ & $0.01932$ & $0.019318838080$ & $0.52$ \\ 
$0.4$ & $0.007099$ & $0.00718$ & $0.007184856672$ & $0.50$%
\end{tabular}
\\ \hline
\end{tabular}%
\end{equation*}%
\begin{equation*}
\begin{tabular}{|l|}
\hline
Table 4 b: Maximum values of $V$ at different time levels for $%
k_{2}=k_{3}=10 $ \\ \hline
\begin{tabular}{ccccc}
$t$ & Present & \cite{2011iki} & \cite{2014uc} & at point \\ 
$0.1$ & $0.144178$ & $0.14306$ & $0.143141957500$ & $0.66$ \\ 
$0.2$ & $0.049030$ & $0.04697$ & $0.047006446750$ & $0.56$ \\ 
$0.3$ & $0.018049$ & $0.01725$ & $0.017260356430$ & $0.52$ \\ 
$0.4$ & $0.006711$ & $0.00641$ & $0.006416614856$ & $0.50$%
\end{tabular}
\\ \hline
\end{tabular}%
\end{equation*}%
\begin{equation*}
\begin{tabular}{|l|}
\hline
Table 4 c: Maximum values of $U$ at different time levels for $%
k_{2}=k_{3}=100$ \\ \hline
\begin{tabular}{ccccc}
$t$ & Present & \cite{2011iki} & \cite{2014uc} & at point \\ 
$0.1$ & $0.039322$ & $0.04175$ & $0.041682987260$ & $0.46$ \\ 
$0.2$ & $0.013495$ & $0.01479$ & $0.014770415340$ & $0.58$ \\ 
$0.3$ & $0.004874$ & $0.00534$ & $0.005337325631$ & $0.54$ \\ 
$0.4$ & $0.001808$ & $0.00198$ & $0.001978065014$ & $0.52$%
\end{tabular}
\\ \hline
\end{tabular}%
\end{equation*}%
\begin{equation*}
\begin{tabular}{|l|}
\hline
Table 4 d: Maximum values of $V$ at different time levels for $%
k_{2}=k_{3}=100$ \\ \hline
\begin{tabular}{ccccc}
$t$ & Present & \cite{2011iki} & \cite{2014uc} & at point \\ 
$0.1$ & $0.053927$ & $0.05065$ & $0.050737669860$ & $0.76$ \\ 
$0.2$ & $0.011531$ & $0.01033$ & $0.010356602970$ & $0.64$ \\ 
$0.3$ & $0.003970$ & $0.00350$ & $0.003517189432$ & $0.56$ \\ 
$0.4$ & $0.001464$ & $0.00129$ & $0.001294450199$ & $0.52$%
\end{tabular}
\\ \hline
\end{tabular}%
\end{equation*}

Figs. 6, 7 and 8 show the numerical results obtained for different time
levels $t\in \lbrack 0,1]$ at $k_{2}=k_{3}=10$ for $U$ and $V$ with
different values of $k_{1}$. From the Figs. 7-9, it can be easily seen that
the numerical solutions $U^{n}$ and $V^{n}$ decay to zero as $t$ and $k_{1}$
increased.%
\begin{equation*}
\begin{array}{c}
\begin{array}{cc}
\FRAME{itbpF}{3.712in}{3.712in}{0in}{}{}{6a.bmp}{\special{language
"Scientific Word";type "GRAPHIC";maintain-aspect-ratio TRUE;display
"USEDEF";valid_file "F";width 3.712in;height 3.712in;depth
0in;original-width 3.7537in;original-height 3.7537in;cropleft "0";croptop
"1";cropright "1";cropbottom "0";filename '6a.bmp';file-properties "XNPEU";}}
& \FRAME{itbpF}{3.712in}{3.712in}{0in}{}{}{6b.bmp}{\special{language
"Scientific Word";type "GRAPHIC";maintain-aspect-ratio TRUE;display
"USEDEF";valid_file "F";width 3.712in;height 3.712in;depth
0in;original-width 3.7537in;original-height 3.7537in;cropleft "0";croptop
"1";cropright "1";cropbottom "0";filename '6b.bmp';file-properties "XNPEU";}}%
\end{array}
\\ 
\text{Fig 6: Num. Sol. }U(x,t)\text{ and }V(x,t)\text{ of Problem 3 at
different time levels for }k_{2}=k_{3}=10\text{ while }k_{1}=1%
\end{array}%
\end{equation*}%
\begin{equation*}
\begin{array}{c}
\begin{array}{cc}
\FRAME{itbpF}{3.712in}{3.712in}{0in}{}{}{7a.bmp}{\special{language
"Scientific Word";type "GRAPHIC";maintain-aspect-ratio TRUE;display
"USEDEF";valid_file "F";width 3.712in;height 3.712in;depth
0in;original-width 3.7537in;original-height 3.7537in;cropleft "0";croptop
"1";cropright "1";cropbottom "0";filename '7a.bmp';file-properties "XNPEU";}}
& \FRAME{itbpF}{3.712in}{3.712in}{0in}{}{}{7b.bmp}{\special{language
"Scientific Word";type "GRAPHIC";maintain-aspect-ratio TRUE;display
"USEDEF";valid_file "F";width 3.712in;height 3.712in;depth
0in;original-width 3.7537in;original-height 3.7537in;cropleft "0";croptop
"1";cropright "1";cropbottom "0";filename '7b.bmp';file-properties "XNPEU";}}%
\end{array}
\\ 
\text{Fig 7: Num. Sol. }U(x,t)\text{ and }V(x,t)\text{ of Problem 3 at
different time levels for }k_{2}=k_{3}=10\text{ while }k_{1}=50%
\end{array}%
\end{equation*}%
\begin{equation*}
\begin{array}{c}
\begin{array}{cc}
\FRAME{itbpF}{3.712in}{3.712in}{0in}{}{}{8a.bmp}{\special{language
"Scientific Word";type "GRAPHIC";maintain-aspect-ratio TRUE;display
"USEDEF";valid_file "F";width 3.712in;height 3.712in;depth
0in;original-width 3.7537in;original-height 3.7537in;cropleft "0";croptop
"1";cropright "1";cropbottom "0";filename '8a.bmp';file-properties "XNPEU";}}
& \FRAME{itbpF}{3.712in}{3.712in}{0in}{}{}{8b.bmp}{\special{language
"Scientific Word";type "GRAPHIC";maintain-aspect-ratio TRUE;display
"USEDEF";valid_file "F";width 3.712in;height 3.712in;depth
0in;original-width 3.7537in;original-height 3.7537in;cropleft "0";croptop
"1";cropright "1";cropbottom "0";filename '8b.bmp';file-properties "XNPEU";}}%
\end{array}
\\ 
\text{Fig 8: Num. Sol. }U(x,t)\text{ and }V(x,t)\text{ of Problem 3 at
different time levels for }k_{2}=k_{3}=10\text{ while }k_{1}=100%
\end{array}%
\end{equation*}

\section{Conclusion}

\qquad Evaluation of the suggested algorithm has done by studying three test
problem. \ The first two of them have analytical solution. For the problem
1, the B-spline($\lambda =0$) collocation algorithm provides the same error
with the suggested algorithm whereas, for the problem 2, errors of the
suggested algorithm is smaller than that of the B-spline collocation method.
The cubic B-spline with additional term defines the extended B-spline
functions. The extended B-spline collocation method with a free parameters
gives smaller error than both the cubic B-spline and cubic trigonometric
cubic B-spline collocation methods. \ Shock propagation is studied \  \ in
the text problem 3.\ Tabulated results and graphical solutions cohere with
results of both the Galerkin quadratic B-spline finite element method and
the modified B-spline collocation method. \ As a result, the trigonometric
B-spline based numerical algorithm can be used for solving the system of
partial differential equation reliably.

\end{document}